\documentclass[leqno,12pt]{amsart}
\usepackage{amssymb} 
\theoremstyle{plain} 
\newtheorem{theorem}{\indent\sc Theorem}[section] 

\newtheorem{proposition}[theorem]{\indent\sc Proposition}

\theoremstyle{definition} 

\begin{document}

\title[Coupled Painlev\'e systems]{Coupled Painlev\'e systems with affine Weyl group symmetry of types $A_7^{(2)},A_5^{(2)}$ and $D_4^{(3)}$ \\}

\renewcommand{\thefootnote}{\fnsymbol{footnote}}
\footnote[0]{2000\textit{ Mathematics Subjet Classification}.
34M55; 34M45; 58F05; 32S65.}

\keywords{ 
Affine Weyl group, birational symmetry, coupled Painlev\'e system.}
\maketitle

\begin{abstract}
 We find a four-parameter family of coupled Painlev\'e VI systems in dimension four with affine Weyl group symmetry of type $A_7^{(2)}$. This is the first example which gave higher-order Painlev\'e equations of type $A_{2l+5}^{(2)}$. We then give an explicit description of a confluence process from this system to a 3-parameter family of coupled Painlev\'e V and III systems in dimension four with $W(A_5^{(2)})$-symmetry. For a degenerate system of $A_5^{(2)}$ system, we also find a two-parameter family of ordinary differential systems in dimension four with affine Weyl group symmetry of type $D_4^{(3)}$.  This is the first example which gave higher-order Painlev\'e equations of type $D_4^{(3)}$. We show that for each system, we give its symmetry and holomorphy conditions. These symmetries, holomorphy conditions and invariant divisors are new.
\end{abstract}

\section{Introduction}
In \cite{Sasa8}, we presented a 4-parameter family of coupled Painlev\'e VI systems in dimension four with affine Weyl group symmetry of type $E_6^{(2)}$.

Now, we consider the following problem.

{\bf Problem} \quad For each affine root system $X_i^{(2)}$ with affine Weyl group $W(X_i^{(2)})$, find a system of differential equations for which $W(X_i^{(2)})$ acts as its B{\"a}cklund transformations.

\vspace{0.1cm}

We will complete the study of the above problem in a series of papers, for which this paper is the second, resulting in a series of equations for the remaining affine root systems of types $A_l^{(2)}$ and $D_4^{(3)}$. This paper is the stage in this project where we find a 4-parameter family of coupled Painlev\'e VI systems in dimension four with $W(A_7^{(2)})$-symmetry explicitly given by
\begin{align}
\begin{split}
\frac{dx}{dt}&=\frac{\partial H}{\partial y}, \quad \frac{dy}{dt}=-\frac{\partial H}{\partial x}, \quad \frac{dz}{dt}=\frac{\partial H}{\partial w}, \quad \frac{dw}{dt}=-\frac{\partial H}{\partial z}
\end{split}
\end{align}
with the polynomial Hamiltonian
\begin{align}
\begin{split}
t(t-1)H =&x^2y^3-((t+1)x+2\alpha_2+\alpha_3+\alpha_4)xy^2\\
&+\{tx^2+(\alpha_0+2\alpha_2+\alpha_3+\alpha_4-(\alpha_0+\alpha_1)t)x\\
&+\alpha_2(\alpha_2+\alpha_3+\alpha_4)\}y+\alpha_1 t x\\
&+\frac{1}{4}[-z^2w^4+2\alpha_4 zw^3+((2t-1)z^2-\alpha_4^2)w^2\\
&-2\{-\alpha_0-2\alpha_2-\alpha_3-\alpha_4+(\alpha_4+1)t\}zw-t(t-1)z^2]\\
&+y((y-1)x-\alpha_2)zw.
\end{split}
\end{align}
Here $x,y,z$ and $w$ denote unknown complex variables, and $\alpha_0,\alpha_1, \dots ,\alpha_4$ are complex parameters satisfying the relation:
\begin{equation}
\alpha_0+\alpha_1+2\alpha_2+2\alpha_3+\alpha_4=1.
\end{equation}
This is the first example which gave higher-order Painlev\'e equations of type $A_{2l+5}^{(2)}$.

We then give an explicit description of a confluence process from this system to a 3-parameter family of ordinary differential systems in dimension four with $W(C_3^{(1)})$-symmetry explicitly given by
\begin{align}
\begin{split}
\frac{dx}{dt}&=\frac{\partial H}{\partial y}, \quad \frac{dy}{dt}=-\frac{\partial H}{\partial x}, \quad \frac{dz}{dt}=\frac{\partial H}{\partial w}, \quad \frac{dw}{dt}=-\frac{\partial H}{\partial z}
\end{split}
\end{align}
with the polynomial Hamiltonian
\begin{align}
\begin{split}
t^2H =&x^2y^3-(tx+2\alpha_1+\alpha_2+\alpha_3)xy^2-\{(\alpha_0 t-1)x-\alpha_1(\alpha_1+\alpha_2+\alpha_3)\}y-tx\\
&-\frac{z^2w^4}{4}+\frac{\alpha_3}{2}zw^3+\frac{1}{4}(2tz^2-\alpha_3^2)w^2-\frac{1}{2}\{(\alpha_3+1)t-1\}zw-\frac{t^2 z^2}{4}\\
&+(xy-\alpha_1)yzw.
\end{split}
\end{align}
Here $x,y,z$ and $w$ denote unknown complex variables, and $\alpha_0,\alpha_1,\alpha_2,\alpha_3$ are complex parameters satisfying the relation:
\begin{equation}
\alpha_0+2\alpha_1+2\alpha_2+\alpha_3=1.
\end{equation}

Moreover, the B{\"a}cklund transformation group for $C_3^{(1)}$ root system can be obtained from that for type $A_7^{(2)}$ by this degeneration process.

We show that by making birational and symplectic transformations with some parameter's change this system is equivalent to a 3-parameter family of coupled Painlev\'e V and III systems in dimension four with $W(A_5^{(2)})$-symmetry explicitly given by
\begin{align}
\begin{split}
\frac{dx}{dt}&=\frac{\partial H}{\partial y}, \quad \frac{dy}{dt}=-\frac{\partial H}{\partial x}, \quad \frac{dz}{dt}=\frac{\partial H}{\partial w}, \quad \frac{dw}{dt}=-\frac{\partial H}{\partial z}
\end{split}
\end{align}
with the polynomial Hamiltonian
\begin{align}
\begin{split}
t^2H =&-tx(x-1)y^2+\{x^2+((\alpha_0+\alpha_1)t-1)x-\alpha_1 t\}y+(\alpha_2+\alpha_3)x-(\alpha_0 t-1)\alpha_1\\
&-\frac{z^2w^4}{4}+\frac{\alpha_3}{2}zw^3+\frac{1}{4}(2tz^2-\alpha_3^2)w^2-\frac{1}{2}\{(\alpha_0+\alpha_1+2\alpha_2+2\alpha_3)t-1\}zw-\frac{t^2 z^2}{4}\\
&-xzw.
\end{split}
\end{align}
Here $x,y,z$ and $w$ denote unknown complex variables, and $\alpha_0,\alpha_1,\alpha_2,\alpha_3$ are complex parameters satisfying the relation:
\begin{equation}
\alpha_0+\alpha_1+2\alpha_2+\alpha_3=1.
\end{equation}
This is the second example which gave higher-order Painlev\'e equations of type $A_{2l+3}^{(2)}$.

For a degenerate system of $A_5^{(2)}$ system, we find a 2-parameter family of ordinary differential systems in dimension four with $W(D_4^{(3)})$-symmetry explicitly given by
\begin{align}
\begin{split}
\frac{dx}{dt}&=\frac{\partial H}{\partial y}, \quad \frac{dy}{dt}=-\frac{\partial H}{\partial x}, \quad \frac{dz}{dt}=\frac{\partial H}{\partial w}, \quad \frac{dw}{dt}=-\frac{\partial H}{\partial z}
\end{split}
\end{align}
with the polynomial Hamiltonian
\begin{align}
\begin{split}
&(4t+3)(16t^2-12t+9)H=\\
&-12tx^2y^4+12(x+2\alpha_0 t)xy^3+3\{(\alpha_0+14\alpha_1+9\alpha_2)x-4\alpha_0^2 t\}y^2\\
&+\{-24tx^2-32t^2 x-3\alpha_0(5\alpha_0+14\alpha_1+9\alpha_2)\}y+24x\{x+(4\alpha_0+7\alpha_1+4\alpha_2)t\}\\
&-12tz^2w^4-12\{z-2(\alpha_0+\alpha_2)t\}zw^3+3\{8t^2 z^2-(\alpha_0+6\alpha_1-3\alpha_2)z-4(\alpha_0+\alpha_2)^2 t\}w^2\\
&+\{12tz^2-8(7\alpha_0+8\alpha_1+7\alpha_2)t^2 z+3(5-6\alpha_1-8\alpha_1^2-4\alpha_2)\}w\\
&-3\{(4t^3+1)z+(19\alpha_0+10\alpha_1+7\alpha_2)t\}z\\
&-12x\{-4txwy^3+(4xw+4txw^2-2tzw^2+zw+2t^2 z)y^2\\
&+(-2xw^2-7tz-4t^2 zw+zw^2+4tzw^3)y+(8tw+5)z\}+24\alpha_0^2 tyw\\
&+12\alpha_1 xw(2w-7y)+6\alpha_2 xyw\{4t(2w-y)-3\}\\
&+6\alpha_0 y\{4t(2x-z)w^2+(x-12t xy+2z)w+4t^2 z\}+24\alpha_0 \alpha_2 tyw.
\end{split}
\end{align}
Here $x,y,z$ and $w$ denote unknown complex variables, and $\alpha_0,\alpha_1,\alpha_2$ are complex parameters satisfying the relation:
\begin{equation}
\alpha_0+2\alpha_1+\alpha_2=1.
\end{equation}
This is the first example which gave higher order Painlev\'e type systems of type $D_{4}^{(3)}$.

We remark that for each system we tried to seek its first integrals of polynomial type with respect to $x,y,z,w$. However, we can not find. Of course, each Hamiltonian is not its first integral.

The B{\"a}cklund transformations of each system (except for $C_3^{(1)}$-system) satisfy
\begin{equation}
s_i(g)=g+\frac{\alpha_i}{f_i}\{f_i,g\}+\frac{1}{2!} \left(\frac{\alpha_i}{f_i} \right)^2 \{f_i,\{f_i,g\} \}+\cdots \quad (g \in {\Bbb C}(t)[x,y,z,w]),
\end{equation}
where poisson bracket $\{,\}$ satisfies the relations:
$$
\{y,x\}=\{w,z\}=1, \quad the \ others \ are \ 0.
$$
Since these B{\"a}cklund transformations have Lie theoretic origin, similarity reduction of a Drinfeld-Sokolov hierarchy admits such a B{\"a}cklund symmetry.

For each system, we also give its holomorphy conditions.

These symmetries, holomorphy conditions and invariant divisors are new.

\section{$A_7^{(2)}$ system}
In this section, we study a 4-parameter family of coupled Painlev\'e VI systems in dimension four with affine Weyl group symmetry of type $A_7^{(2)}$ given by
\begin{align}\label{11}
\begin{split}
\frac{dx}{dt}&=\frac{\partial H}{\partial y}, \quad \frac{dy}{dt}=-\frac{\partial H}{\partial x}, \quad \frac{dz}{dt}=\frac{\partial H}{\partial w}, \quad \frac{dw}{dt}=-\frac{\partial H}{\partial z}
\end{split}
\end{align}
with the polynomial Hamiltonian
\begin{align}\label{12}
\begin{split}
t(t-1)H =&x^2y^3-((t+1)x+2\alpha_2+\alpha_3+\alpha_4)xy^2\\
&+\{tx^2+(\alpha_0+2\alpha_2+\alpha_3+\alpha_4-(\alpha_0+\alpha_1)t)x\\
&+\alpha_2(\alpha_2+\alpha_3+\alpha_4)\}y+\alpha_1 t x\\
&+\frac{1}{4}[-z^2w^4+2\alpha_4 zw^3+((2t-1)z^2-\alpha_4^2)w^2\\
&-2\{-\alpha_0-2\alpha_2-\alpha_3-\alpha_4+(\alpha_4+1)t\}zw-t(t-1)z^2]\\
&+y((y-1)x-\alpha_2)zw.
\end{split}
\end{align}
Here $x,y,z$ and $w$ denote unknown complex variables, and $\alpha_0,\alpha_1, \dots ,\alpha_4$ are complex parameters satisfying the relation:
\begin{equation}\label{13}
\alpha_0+\alpha_1+2\alpha_2+2\alpha_3+\alpha_4=1.
\end{equation}
This is the first example which gave higher order Painlev\'e type systems of type $A_{7}^{(2)}$.

We remark that for this system we tried to seek its first integrals of polynomial type with respect to $x,y,z,w$. However, we can not find. Of course, the Hamiltonian $H$ is not the first integral.

It is known that the Painlev\'e VI system admits the affine Weyl group symmetry of type $B_4^{(1)}$ as the group of its B{\"a}cklund transformations in addition to the diagram automorphisms of type $D_4^{(1)}$. The diagram automorphisms change the time variable $t$. However, the system \eqref{11} admits the affine Weyl group symmetry of type $A_7^{(2)}$ as the group of its B{\"a}cklund transformations, whose generators $s_0,s_1,\ldots,s_4$ are determined by the invariant divisors \eqref{19}. Of course, these transformations do not change the time variable $t$.

We also remark that the iniariant divisors of the system \eqref{11} are different from the ones of a 4-parameter family of 2-coupled Painlev\'e VI systems in dimension four with affine Weyl group symmetry of type $E_6^{(2)}$ given in the paper \cite{Sasa8}.

Now, we show that each principal part of this Hamiltonian can be transformed into canonical Painlev\'e VI Hamiltonian by birational and symplectic transformations.

At first, we study the Hamiltonian system
\begin{align}
\begin{split}
\frac{dx}{dt}&=\frac{\partial K_1}{\partial y}, \quad \frac{dy}{dt}=-\frac{\partial K_1}{\partial x}
\end{split}\label{14}
\end{align}
with the polynomial Hamiltonian
\begin{align}
\begin{split}
t(1-t)K_1 =&x^2y^3-((t+1)x+2\alpha_2+\alpha_3+\alpha_4)xy^2\\
&+\{tx^2+(\alpha_0+2\alpha_2+\alpha_3+\alpha_4-(\alpha_0+\alpha_1)t)x\\
&+\alpha_2(\alpha_2+\alpha_3+\alpha_4)\}y+\alpha_1 t x,
\end{split}
\end{align}
where setting $z=w=0$ in the Hamiltonian $H$, we obtain $K_1$.

We transform the Hamiltonian \eqref{14} into the Painlev\'e VI Hamiltonian:
\begin{align}\label{HPVI}
\begin{split}
&H_{VI}(x,y,t;\beta_0,\beta_1,\beta_2,\beta_3,\beta_4)\\
&=\frac{1}{t(t-1)}[y^2(x-t)(x-1)x-\{(\beta_0-1)(x-1)x+\beta_3(x-t)x\\
&+\beta_4(x-t)(x-1)\}y+\beta_2(\beta_1+\beta_2)x]  \quad (\beta_0+\beta_1+2\beta_2+\beta_3+\beta_4=1). 
\end{split}
\end{align}

{\bf Step 1:} We make the change of variables:
\begin{equation}
x_1=y, \quad y_1=-x.
\end{equation}
Then, we can obtain the Painlev\'e VI Hamiltonian:
\begin{align}
\begin{split}
&H_{VI}(x_1,y_1,t;\alpha_3,\alpha_3+\alpha_4,\alpha_2,\alpha_0,\alpha_1).
\end{split}
\end{align}
Of course, the parameters $\alpha_i$ and $\beta_j$ satisfy the relations:
\begin{equation}
\beta_0+\beta_1+2\beta_2+\beta_3+\beta_4=\alpha_0+\alpha_1+2\alpha_2+2\alpha_3+\alpha_4=1.
\end{equation}
We remark that all transformations are symplectic.

Next, we study the Hamiltonian system
\begin{align}
\begin{split}
\frac{dz}{dt}&=\frac{\partial K_2}{\partial w}, \quad \frac{dw}{dt}=-\frac{\partial K_2}{\partial z}
\end{split}
\end{align}
with the polynomial Hamiltonian
\begin{align}\label{15}
\begin{split}
t(1-t)K_2 =&\frac{1}{4}[-z^2w^4+2\alpha_4 zw^3+((2t-1)z^2-\alpha_4^2)w^2\\
&-2\{-\alpha_0-2\alpha_2-\alpha_3-\alpha_4+(\alpha_4+1)t\}zw-t(t-1)z^2],
\end{split}
\end{align}
where setting $x=y=0$ in the Hamiltonian $H$, we obtain $K_2$.

Let us transform the Hamiltonian \eqref{15} into the Painlev\'e VI Hamiltonian.

{\bf Step 1:} We make the change of variables:
\begin{equation}
z_1=2 \sqrt{t-1}z, \quad w_1=\frac{1}{2 \sqrt{t-1}}w+\frac{1}{2}.
\end{equation}
By this transformation, in the coordinate system $(Z_1,W_1)=(1/z_1,w_1)$  two of four accessible singular points are transformed into $W_1=0$ and $W_1=1$.

{\bf Step 2:} We make the change of variables:
\begin{equation}
z_2=-(z_1w_1-\alpha_4)w_1, \quad w_2=\frac{1}{w_1}.
\end{equation}

{\bf Step 3:} We make the change of variables:
\begin{align}
\begin{split}
&z_3=\frac{z_2}{2t-1+2\sqrt{t(t-1)}}, \quad w_3=(2t-1+2\sqrt{t(t-1)})w_2+2(1-t-\sqrt{t(t-1)}),\\
&t=\frac{2T(T-1)-(2T-1)\sqrt{T(T-1)}}{4T(T-1)}.
\end{split}
\end{align}
By this transformation, in the coordinate system $(Z_2,W_2)=(1/z_3,w_3)$  the others are transformed into $W_2=0$ and $W_2=\frac{1}{T}$. We remark that it is not $W_2 = \infty$ but  $W_2=0$ because we consider in the coordinate system $(z_2,w_2)$.

{\bf Step 4:} We make the change of variables:
\begin{equation}
z_4=-(z_3w_3-\alpha_4)w_3, \quad w_4=\frac{1}{w_3}.
\end{equation}

{\bf Step 5:} We make the change of variables:
\begin{equation}
z_5=w_4, \quad w_5=-z_4.
\end{equation}
Then, we can obtain the Painlev\'e VI Hamiltonian:
\begin{align}
\begin{split}
&-\frac{1}{2} H_{VI}(z_5,w_5,T;\alpha_1+\alpha_3-1,\alpha_1+\alpha_3,\alpha_4,\alpha_0+2\alpha_2+\alpha_3,\alpha_0+2\alpha_2+\alpha_3).
\end{split}
\end{align}
Of course, the parameters $\alpha_i$ and $\beta_j$ satisfy the relations:
\begin{equation}
\beta_0+\beta_1+2\beta_2+\beta_3+\beta_4=2(\alpha_0+\alpha_1+2\alpha_2+2\alpha_3+\alpha_4)-1=1.
\end{equation}
We remark that all transformations are symplectic.

\begin{figure}
\unitlength 0.1in
\begin{picture}(43.30,21.70)(24.05,-23.35)
%
\special{pn 20}%
\special{ar 2830 590 425 425  0.0000000 6.2831853}%
%
\special{pn 20}%
\special{ar 2840 1910 425 425  0.0000000 6.2831853}%
%
\special{pn 20}%
\special{ar 3880 1260 425 425  0.0000000 6.2831853}%
%
\special{pn 20}%
\special{ar 5120 1260 425 425  0.0000000 6.2831853}%
%
\special{pn 20}%
\special{ar 6310 1260 425 425  0.0000000 6.2831853}%
%
\special{pn 20}%
\special{pa 3190 820}%
\special{pa 3510 990}%
\special{fp}%
%
\special{pn 20}%
\special{pa 3190 1670}%
\special{pa 3530 1500}%
\special{fp}%
%
\special{pn 20}%
\special{pa 4300 1250}%
\special{pa 4670 1250}%
\special{fp}%
%
\special{pn 20}%
\special{pa 5890 1160}%
\special{pa 5560 1160}%
\special{fp}%
\special{sh 1}%
\special{pa 5560 1160}%
\special{pa 5627 1180}%
\special{pa 5613 1160}%
\special{pa 5627 1140}%
\special{pa 5560 1160}%
\special{fp}%
%
\special{pn 20}%
\special{pa 5900 1410}%
\special{pa 5550 1410}%
\special{fp}%
\special{sh 1}%
\special{pa 5550 1410}%
\special{pa 5617 1430}%
\special{pa 5603 1410}%
\special{pa 5617 1390}%
\special{pa 5550 1410}%
\special{fp}%
\put(26.1000,-20.0000){\makebox(0,0)[lb]{$y-1$}}%
\put(26.6000,-6.8000){\makebox(0,0)[lb]{$y$}}%
\put(37.4000,-13.5000){\makebox(0,0)[lb]{$x$}}%
\put(61.7000,-13.4000){\makebox(0,0)[lb]{$z$}}%
\put(47.3000,-13.6000){\makebox(0,0)[lb]{$y+w^2-t$}}%
\end{picture}%
\label{fig:A72}
\caption{This figure denotes Dynkin diagram of type $A_7^{(2)}$. The symbol in each circle denotes the invariant divisors of the system \eqref{11} (see Theorem \ref{th:1.1}).}
\end{figure}
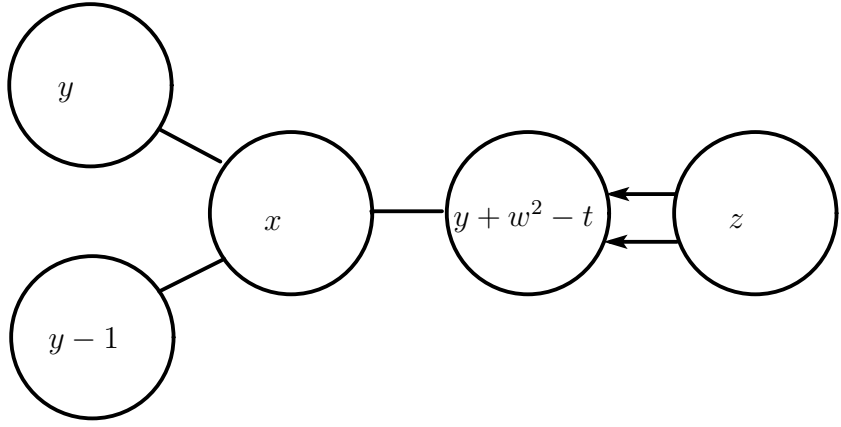

\begin{theorem}\label{th:1.1}
The system \eqref{11} admits extended affine Weyl group symmetry of type $A_7^{(2)}$ as the group of its B{\"a}cklund transformations, whose generators $s_0,s_1,\ldots,s_4,\pi$ defined as follows$:$ with {\it the notation} $(*):=(x,y,z,w,t;\alpha_0,\alpha_1,\ldots,\alpha_4)$\rm{: \rm}
\begin{align}
\begin{split}
s_0:(*) \rightarrow &\left(x+\frac{\alpha_0}{y-1},y,z,w,t;-\alpha_0,\alpha_1,\alpha_2+\alpha_0,\alpha_3,\alpha_4 \right),\\
s_1:(*) \rightarrow &\left(x+\frac{\alpha_1}{y},y,z,w,t;\alpha_0,-\alpha_1,\alpha_2+\alpha_1,\alpha_3,\alpha_4 \right),\\
s_2:(*) \rightarrow &\left(x,y-\frac{\alpha_2}{x},z,w,t;\alpha_0+\alpha_2,\alpha_1+\alpha_2,-\alpha_2,\alpha_3+\alpha_2,\alpha_4 \right),\\
s_3:(*) \rightarrow &\left(x+\frac{\alpha_3}{y+w^2-t},y,z+\frac{2\alpha_3 w}{y+w^2-t},w,t;\alpha_0,\alpha_1,\alpha_2+\alpha_3,-\alpha_3,\alpha_4+2\alpha_3 \right),\\
s_4:(*) \rightarrow &\left(x,y,z,w-\frac{\alpha_4}{z},t;\alpha_0,\alpha_1,\alpha_2,\alpha_3+\alpha_4,-\alpha_4 \right),\\
\pi:(*) \rightarrow &\left(-x,1-y,\sqrt{-1}z,\frac{w}{\sqrt{-1}},1-t;\alpha_1,\alpha_0,\alpha_2,\alpha_3,\alpha_4 \right).
\end{split}
\end{align}
\end{theorem}
We note that the B{\"a}cklund transformations of this system satisfy
\begin{equation}\label{universal}
s_i(g)=g+\frac{\alpha_i}{f_i}\{f_i,g\}+\frac{1}{2!} \left(\frac{\alpha_i}{f_i} \right)^2 \{f_i,\{f_i,g\} \}+\cdots \quad (g \in {\Bbb C}(t)[x,y,z,w]),
\end{equation}
where poisson bracket $\{,\}$ satisfies the relations:
$$
\{y,x\}=\{w,z\}=1, \quad the \ others \ are \ 0.
$$
Since these B{\"a}cklund transformations have Lie theoretic origin, similarity reduction of a Drinfeld-Sokolov hierarchy admits such a B{\"a}cklund symmetry.

\begin{proposition}
This system has the following invariant divisors\rm{:\rm}
\begin{center}\label{19}
\begin{tabular}{|c|c|c|} \hline
parameter's relation & $f_i$ \\ \hline
$\alpha_0=0$ & $f_0:=y-1$  \\ \hline
$\alpha_1=0$ & $f_1:=y$  \\ \hline
$\alpha_2=0$ & $f_2:=x$  \\ \hline
$\alpha_3=0$ & $f_3:=y+w^2-t$  \\ \hline
$\alpha_4=0$ & $f_4:=z$  \\ \hline
\end{tabular}
\end{center}
\end{proposition}
We note that when $\alpha_1=0$, we see that the system \eqref{11} admits a particular solution $y=0$, and when $\alpha_3=0$, after we make the birational and symplectic transformations:
\begin{equation}
x_3=x, \ y_3=y+w^2-t, \ z_3=z-2xw, \ w_3=w
\end{equation}
we see that the system \eqref{11} admits a particular solution $y_3=0$.

\begin{proposition}
Let us define the following translation operators{\rm : \rm}
\begin{align}
\begin{split}
&T_1:=\pi s_0 s_2 s_3 s_4 s_3 s_2 s_0, \quad T_2:=s_0 T_1 s_0, \quad T_3:=s_2 T_2 s_2, \quad T_4:=s_4 s_3 T_3 s_3 s_4.
\end{split}
\end{align}
These translation operators act on parameters $\alpha_i$ as follows$:$
\begin{align}
\begin{split}
T_1(\alpha_0,\alpha_1,\ldots,\alpha_4)=&(\alpha_0,\alpha_1,\ldots,\alpha_4)+(-1,1,0,0,0),\\
T_2(\alpha_0,\alpha_1,\ldots,\alpha_4)=&(\alpha_0,\alpha_1,\ldots,\alpha_4)+(1,1,-1,0,0),\\
T_3(\alpha_0,\alpha_1,\ldots,\alpha_4)=&(\alpha_0,\alpha_1,\ldots,\alpha_4)+(0,0,1,-1,0),\\
T_4(\alpha_0,\alpha_1,\ldots,\alpha_4)=&(\alpha_0,\alpha_1,\ldots,\alpha_4)+(0,0,0,-1,2).
\end{split}
\end{align}
\end{proposition}

\begin{theorem}\label{th:1.2}
Let us consider a polynomial Hamiltonian system with Hamiltonian $K \in {\Bbb C}(t)[x,y,z,w]$. We assume that

$(A1)$ $deg(K)=6$ with respect to $x,y,z,w$.

$(A2)$ This system becomes again a polynomial Hamiltonian system in each coordinate system $r_i \ (i=0,1,\ldots,4)${\rm : \rm}
\begin{align}
\begin{split}
r_0:&x_0=\frac{1}{x}, \ y_0=-((y-1)x+\alpha_0)x, \ z_0=z, \ w_0=w, \\
r_1:&x_1=\frac{1}{x}, \ y_1=-(yx+\alpha_1)x, \ z_1=z, \ w_1=w, \\
r_2:&x_2=-(xy-\alpha_2)y, \ y_2=\frac{1}{y}, \ z_2=z, \ w_2=w,\\
r_3:&x_3=\frac{1}{x}, \ y_3=-\left((y+w^2-t)x+\alpha_3 \right)x, \ z_3=z-2xw, \ w_3=w, \\
r_4:&x_4=x, \ y_4=y, \ z_4=-(zw-\alpha_4)w, \ w_4=\frac{1}{w}
\end{split}
\end{align}
Then such a system coincides with the system \eqref{11} with the polynomial Hamiltonian \eqref{12}.
\end{theorem}
By this theorem, we can also recover the parameter's relation \eqref{13}.

We note that the condition $(A2)$ should be read that
\begin{align*}
&r_j(K) \quad (j=0,1,2,4), \quad r_3(K+x)
\end{align*}
are polynomials with respect to $x,y,z,w$.

\section{$C_3^{(1)}$ system}

In this section, we study a 3-parameter family of ordinary differential systems in dimension four with affine Weyl group symmetry of type $C_3^{(1)}$. At first, we consider an explicit description of a confluence process from the system \eqref{11}.

\begin{theorem}\label{6.0}
For the system \eqref{11} of type $A_7^{(2)}$, we make the change of parameters and variables
\begin{gather}
\begin{gathered}\label{65}
\alpha_0=\frac{1}{\varepsilon}+A_0, \quad \alpha_1=-\frac{1}{\varepsilon}, \quad \alpha_2=A_1, \quad \alpha_3=A_2, \quad \alpha_4=A_3,\\
\end{gathered}\\
\begin{gathered}\label{66}
t=\frac{T}{\varepsilon}, \quad x=\varepsilon X,\quad y=\frac{Y}{\varepsilon}, \quad z=\sqrt{\varepsilon} Z, \quad w=\frac{W}{\sqrt{\varepsilon}}
\end{gathered}
\end{gather}
from $\alpha_0,\alpha_1,\alpha_2,\alpha_3,x,y,z,w$ to $A_0,A_1,A_2,A_3,X,Y,Z,W$. Then the system \eqref{11} can also be written in the new variables $X,Y,Z,W$ and parameters $A_0,A_1,A_2,A_3$ as a Hamiltonian system. This new system tends to
\begin{align}\label{61}
\begin{split}
\frac{dx}{dt}&=\frac{\partial H}{\partial y}, \quad \frac{dy}{dt}=-\frac{\partial H}{\partial x}, \quad \frac{dz}{dt}=\frac{\partial H}{\partial w}, \quad \frac{dw}{dt}=-\frac{\partial H}{\partial z}
\end{split}
\end{align}
with the polynomial Hamiltonian
\begin{align}\label{62}
\begin{split}
t^2H =&x^2y^3-(tx+2\alpha_1+\alpha_2+\alpha_3)xy^2-\{(\alpha_0 t-1)x-\alpha_1(\alpha_1+\alpha_2+\alpha_3)\}y-tx\\
&-\frac{z^2w^4}{4}+\frac{\alpha_3}{2}zw^3+\frac{1}{4}(2tz^2-\alpha_3^2)w^2-\frac{1}{2}\{(\alpha_3+1)t-1\}zw-\frac{t^2 z^2}{4}\\
&+(xy-\alpha_1)yzw,
\end{split}
\end{align}
where, for notational convenience, we have renamed $A_i,X,Y,Z,W$ to $\alpha_i,x,y,z,w$ (which are not the same as the previous $\alpha_i,x,y,z,w$).
\end{theorem}
Here $x,y,z$ and $w$ denote unknown complex variables, and $\alpha_0,\alpha_1,\alpha_2,\alpha_3$ are complex parameters satisfying the relation:
\begin{equation}\label{63}
\alpha_0+2\alpha_1+2\alpha_2+\alpha_3=1.
\end{equation}
We remark that for this system we tried to seek its first integrals of polynomial type with respect to $x,y,z,w$. However, we can not find. Of course, the Hamiltonian $H$ is not the first integral.

\begin{theorem}\label{th:6.1}
The system \eqref{61} admits the affine Weyl group symmetry of type $C_3^{(1)}$ as the group of its B{\"a}cklund transformations, whose generators $S_0,S_1,S_2,S_3$ defined as follows$:$ with {\it the notation} $(*):=(x,y,z,w,t;\alpha_0,\alpha_1,\alpha_2,\alpha_3)$\rm{: \rm}
\begin{align}
\begin{split}
S_0:(*) \rightarrow &\left(-\left(x+\frac{\alpha_0}{y}+\frac{1}{y^2} \right),-y,\sqrt{-1}z,\frac{w}{\sqrt{-1}},-t;-\alpha_0,\alpha_1+\alpha_0,\alpha_2,\alpha_3 \right),\\
S_1:(*) \rightarrow &\left(x,y-\frac{\alpha_1}{x},z,w,t;\alpha_0+2\alpha_1,-\alpha_1,\alpha_2+\alpha_1,\alpha_3 \right),\\
S_2:(*) \rightarrow &\left(x+\frac{\alpha_2}{y+w^2-t},y,z+\frac{2\alpha_2 w}{y+w^2-t},w,t;\alpha_0,\alpha_1+\alpha_2,-\alpha_2,\alpha_3+2\alpha_2 \right),\\
S_3:(*) \rightarrow &\left(x,y,z,w-\frac{\alpha_3}{z},t;\alpha_0,\alpha_1,\alpha_2+\alpha_3,-\alpha_3 \right).
\end{split}
\end{align}
\end{theorem}
We remark that the transformations $S_i \ (i=1,2,3)$ satisfy the relation \eqref{universal} (see Section 2). However, the transformation $S_0$ does not satisfy so.

\begin{theorem}\label{6.6}
For the transformations \eqref{65}, \eqref{66} given in Theorem \rm{ \ref{6.0}, \rm} we can choose a subgroup $W_{A_7^{(2)} \rightarrow C_3^{(1)}}$ of the B{\"a}cklund transformation group $W(A_7^{(2)})$ so that $W_{A_7^{(2)} \rightarrow C_3^{(1)}}$ converges to the B{\"a}cklund transformation group $W(C_3^{(1)})$ of the system \eqref{61}.
\end{theorem}
{\it Proof.} \quad Notice that
$$
\alpha_0+\alpha_1+2\alpha_2+2\alpha_3+\alpha_4=A_0+2A_1+2A_2+A_3=1
$$
and the change of variables from $(x,y,z,w)$ to $(X,Y,Z,W)$ is symplectic. Choose $S_i \ (i=0,1,2,3)$ as
$$
S_0:=s_0s_1, \ S_1:=s_2, \ S_2:=s_3, \ S_3:=s_4,
$$
where $s_i$ are given by Theorem \ref{th:1.1}. Then the transformations $S_i$ are reflections of the parameters $A_0,A_1,A_2,A_3$. The transformation group $W_{A_7^{(2)} \rightarrow C_3^{(1)}}=<S_0,S_1,S_2,S_3>$ coincides with the group given in Theorem \ref{th:6.1} as $\varepsilon \rightarrow 0$. \qed

\begin{proposition}
This system has the following invariant divisors\rm{:\rm}
\begin{center}
\begin{tabular}{|c|c|c|} \hline
parameter's relation & $f_i$ \\ \hline
$\alpha_1=0$ & $f_1:=x$  \\ \hline
$\alpha_2=0$ & $f_2:=y+w^2-t$  \\ \hline
$\alpha_3=0$ & $f_3:=z$  \\ \hline
\end{tabular}
\end{center}
\end{proposition}

\begin{theorem}\label{th:6.2}
Let us consider a polynomial Hamiltonian system with Hamiltonian $K \in {\Bbb C}(t)[x,y,z,w]$. We assume that

$(G1)$ $deg(K)=6$ with respect to $x,y,z,w$.

$(G2)$ This system becomes again a polynomial Hamiltonian system in each coordinate system $r_i \ (i=0,1,3)${\rm : \rm}
\begin{align}
\begin{split}
r_0:&x_0=x+\frac{\alpha_0}{y}+\frac{1}{y^2}, \ y_0=y, \ z_0=z, \ w_0=w, \\
r_1:&x_1=-(xy-\alpha_1)y, \ y_1=\frac{1}{y}, \ z_1=z, \ w_1=w, \\
r_2:&x_2=\frac{1}{x}, \ y_2=-\left((y+w^2-t)x+\alpha_2 \right)x, \ z_2=z-2xw, \ w_2=w,\\
r_3:&x_3=x, \ y_3=y, \ z_3=-(zw-\alpha_3)w, \ w_3=\frac{1}{w}.
\end{split}
\end{align}
Then such a system coincides with the system \eqref{61} with the polynomial Hamiltonian \eqref{62}.
\end{theorem}
By this theorem, we can also recover the parameter's relation \eqref{63}.

We note that the conditions $(G2)$ should be read that
\begin{align*}
&r_j(K) \quad (j=0,1,3), \quad r_2(K+x)
\end{align*}
are polynomials with respect to $x,y,z,w$.

\begin{figure}[ht]
\unitlength 0.1in
\begin{picture}(32.50,21.65)(24.05,-23.30)
%
\special{pn 20}%
\special{ar 2830 590 425 425  0.0000000 6.2831853}%
%
\special{pn 20}%
\special{ar 2855 1905 425 425  0.0000000 6.2831853}%
%
\special{pn 20}%
\special{ar 3960 1225 500 425  0.0000000 6.2831853}%
%
\special{pn 20}%
\special{ar 5230 1235 425 425  0.0000000 6.2831853}%
%
\special{pn 20}%
\special{pa 3190 820}%
\special{pa 3510 990}%
\special{fp}%
%
\special{pn 20}%
\special{pa 3205 1665}%
\special{pa 3545 1495}%
\special{fp}%
%
\special{pn 20}%
\special{pa 4810 1135}%
\special{pa 4480 1135}%
\special{fp}%
\special{sh 1}%
\special{pa 4480 1135}%
\special{pa 4547 1155}%
\special{pa 4533 1135}%
\special{pa 4547 1115}%
\special{pa 4480 1135}%
\special{fp}%
%
\special{pn 20}%
\special{pa 4820 1385}%
\special{pa 4470 1385}%
\special{fp}%
\special{sh 1}%
\special{pa 4470 1385}%
\special{pa 4537 1405}%
\special{pa 4523 1385}%
\special{pa 4537 1365}%
\special{pa 4470 1385}%
\special{fp}%
\put(26.2500,-19.9500){\makebox(0,0)[lb]{$x-1$}}%
\put(26.6000,-6.8000){\makebox(0,0)[lb]{$x$}}%
\put(50.9000,-13.1500){\makebox(0,0)[lb]{$z$}}%
\put(35.1000,-13.3000){\makebox(0,0)[lb]{$yw^2-ty+1$}}%
\end{picture}%
\label{fig:A62}
\caption{This figure denotes Dynkin diagram of type $A_5^{(2)}$. The symbol in each circle denotes the invariant divisors of the system \eqref{31} (see Theorem \ref{th:3.1}).}
\end{figure}

\section{$A_5^{(2)}$ system}

In this section, we study a 3-parameter family of coupled Painlev\'e V and III systems in dimension four with affine Weyl group symmetry of type $A_5^{(2)}$. At first, by making birational and symplectic transformations with some parameter's change the system \eqref{61} is equivalent to the system of type $A_5^{(2)}$.

\begin{theorem}\label{3.0}
For the system \eqref{61} of type $C_3^{(1)}$, we make the change of parameters and variables
\begin{gather}
\begin{gathered}\label{35}
\alpha_0=A_0-A_1, \quad \alpha_1=A_1, \quad \alpha_2=A_2, \quad \alpha_3=A_3,\\
\end{gathered}\\
\begin{gathered}\label{36}
x=\frac{1}{X},\quad y=-(XY+A_1)X, \quad z=Z, \quad w=W
\end{gathered}
\end{gather}
from $\alpha_0,\alpha_1,\alpha_2,\alpha_3,x,y,z,w$ to $A_0,A_1,A_2,A_3,X,Y,Z,W$. Then the system \eqref{61} can also be written in the new variables $X,Y,Z,W$ and parameters $A_0,A_1,A_2,A_3$ as the polynomial Hamiltonian system given by
\begin{align}\label{31}
\begin{split}
\frac{dx}{dt}&=\frac{\partial H}{\partial y}, \quad \frac{dy}{dt}=-\frac{\partial H}{\partial x}, \quad \frac{dz}{dt}=\frac{\partial H}{\partial w}, \quad \frac{dw}{dt}=-\frac{\partial H}{\partial z}
\end{split}
\end{align}
with the polynomial Hamiltonian
\begin{align}\label{32}
\begin{split}
t^2H =&-tx(x-1)y^2+\{x^2+((\alpha_0+\alpha_1)t-1)x-\alpha_1 t\}y+(\alpha_2+\alpha_3)x-(\alpha_0 t-1)\alpha_1\\
&-\frac{z^2w^4}{4}+\frac{\alpha_3}{2}zw^3+\frac{1}{4}(2tz^2-\alpha_3^2)w^2-\frac{1}{2}\{(\alpha_0+\alpha_1+2\alpha_2+2\alpha_3)t-1\}zw-\frac{t^2 z^2}{4}\\
&-xzw,
\end{split}
\end{align}
where, for notational convenience, we have renamed $A_i,X,Y,Z,W$ to $\alpha_i,x,y,z,w$ (which are not the same as the previous $\alpha_i,x,y,z,w$).
\end{theorem}
Here $x,y,z$ and $w$ denote unknown complex variables, and $\alpha_0,\alpha_1,\alpha_2,\alpha_3$ are complex parameters satisfying the relation:
\begin{equation}\label{33}
\alpha_0+\alpha_1+2\alpha_2+\alpha_3=1.
\end{equation}
This is the second example which gave higher order Painlev\'e type systems of type $A_{5}^{(2)}$ (see \cite{Sasa10}). In this case the iniariant divisors are different from the ones of a 3-parameter family of 2-coupled Painlev\'e III systems in dimension four given in the paper \cite{Sasa10}.
\begin{center}
\begin{tabular}{|c||c|c|c|c|} \hline 
Invariant divisors &  $f_0$ & $f_1$ & $f_2$ & $f_3$ \\ \hline
System \eqref{31} &  $x-1$ & $x$ & $yw^2-ty+1$ & $z$ \\ \hline
2-CPIII &  $y$ & $w$ & $xz+t$ & $y+w-1$ \\ \hline
\end{tabular}
\end{center}

We remark that for this system we tried to seek its first integrals of polynomial type with respect to $x,y,z,w$. However, we can not find. Of course, the Hamiltonian $H$ is not the first integral.

We will show that each principal part of this Hamiltonian can be transformed into canonical Painlev\'e V and III Hamiltonian by birational and symplectic transformations. At first, we study the Hamiltonian system
\begin{align}
\begin{split}
\frac{dx}{dt}&=\frac{\partial K_3}{\partial y}, \quad \frac{dy}{dt}=-\frac{\partial K_3}{\partial x}
\end{split}\label{338}
\end{align}
with the polynomial Hamiltonian
\begin{align}
\begin{split}
t^2K_3 =&-tx(x-1)y^2+\{x^2+((\alpha_0+\alpha_1)t-1)x-\alpha_1 t\}y+(\alpha_2+\alpha_3)x-(\alpha_0 t-1)\alpha_1,
\end{split}
\end{align}
where setting $z=w=0$ in the Hamiltonian $H$, we obtain $K_3$.

We transform the Hamiltonian \eqref{338} into the Painlev\'e V Hamiltonian (see \cite{Sasa9}):
\begin{align}\label{HPV}
\begin{split}
&-tH_{V}(x,y,t;\beta_1,\beta_2,\beta_3)\\
&=-ty+x^3y^2-x^2y^2+(\beta_1+2\beta_2)x^2y+(t-1+2\beta_3)xy+\beta_2(\beta_1+\beta_2)x. 
\end{split}
\end{align}

{\bf Step 1:} We make the change of variables:
\begin{equation}
x_1=\frac{y}{t}, \quad y_1=-tx, \quad T=-\frac{1}{t}.
\end{equation}
Then, we can obtain the Painlev\'e V Hamiltonian:
\begin{align}
\begin{split}
&H_{V}\left(x_1,y_1,T;\alpha_2+\alpha_3,\alpha_1,\frac{\alpha_0}{2} \right).
\end{split}
\end{align}
We remark that all transformations are symplectic.

Next, we study the Hamiltonian system
\begin{align}\label{339}
\begin{split}
\frac{dz}{dt}&=\frac{\partial K_4}{\partial w}, \quad \frac{dw}{dt}=-\frac{\partial K_4}{\partial z}
\end{split}
\end{align}
with the polynomial Hamiltonian
\begin{align}\label{337}
\begin{split}
t^2 K_4=&-\frac{z^2w^4}{4}+\frac{\alpha_3}{2}zw^3+\frac{1}{4}(2tz^2-\alpha_3^2)w^2-\frac{1}{2}\{(\alpha_0+\alpha_1+2\alpha_2+2\alpha_3)t-1\}zw-\frac{t^2 z^2}{4},
\end{split}
\end{align}
where setting $x=y=0$ in the Hamiltonian $H$, we obtain $K_4$.

Let us transform the Hamiltonian \eqref{339} into the Painlev\'e III Hamiltonian:
\begin{align}\label{HPIII}
\begin{split}
&H_{III}(q,p,t;\alpha_1,\beta_1)=\frac{q^2p^2-(q^2-(\alpha_1+\beta_1)q-t)p-\alpha_1 q}{t}. 
\end{split}
\end{align}

{\bf Step 0:} We make the change of variables:
\begin{equation}
t=\frac{1}{T_1}.
\end{equation}
We note that
\begin{equation}
dK_4 \wedge dt=-\frac{1}{{T_1}^2}d{\tilde K}_4 \wedge dT_1.
\end{equation}

{\bf Step 1:} We make the change of variables:
\begin{equation}
z_1=\frac{z}{2\sqrt{T_1}}, \quad w_1=2\sqrt{T_1} w+2.
\end{equation}

{\bf Step 2:} We make the change of variables:
\begin{equation}
z_2=-(z_1w_1-\alpha_3)w_1, \quad w_2=\frac{1}{w_1}.
\end{equation}

{\bf Step 3:} We make the change of variables:
\begin{align}
\begin{split}
&z_3=-\frac{1}{T_1}z_2+1, \quad w_3=-T_1 w_2.
\end{split}
\end{align}

{\bf Step 4:} We make the change of variables:
\begin{equation}
z_4=w_3+\frac{T_1}{4}, \quad w_4=-z_3, \quad T_1=4T_2.
\end{equation}

{\bf Step 5:} We make the change of variables:
\begin{equation}
z_5=z_4, \quad w_5=w_4+1,\quad T_2=\sqrt{T_3}.
\end{equation}
Then, we can obtain the Painlev\'e III Hamiltonian:
\begin{align}
\begin{split}
&\frac{1}{2} H_{III}(z_5,w_5,T_3;0,\alpha_3+1).
\end{split}
\end{align}
We remark that all transformations are symplectic.

\begin{theorem}\label{th:3.1}
The system \eqref{31} admits extended affine Weyl group symmetry of type $A_5^{(2)}$ as the group of its B{\"a}cklund transformations, whose generators $s_0,s_1,s_2,s_3,\pi$ defined as follows$:$ with {\it the notation} $(*):=(x,y,z,w,t;\alpha_0,\alpha_1,\alpha_2,\alpha_3)$\rm{: \rm}
\begin{align}
\begin{split}
s_0:(*) \rightarrow &\left(x,y-\frac{\alpha_0}{x-1},z,w,t;-\alpha_0,\alpha_1,\alpha_2+\alpha_0,\alpha_3 \right),\\
s_1:(*) \rightarrow &\left(x,y-\frac{\alpha_1}{x},z,w,t;\alpha_0,-\alpha_1,\alpha_2+\alpha_1,\alpha_3 \right),\\
s_2:(*) \rightarrow &\left(x+\frac{\alpha_2(w^2-t)}{yw^2-ty+1},y,z+\frac{2\alpha_2 yw}{yw^2-ty+1},w,t;\alpha_0+\alpha_2,\alpha_1+\alpha_2,-\alpha_2,\alpha_3+2\alpha_2 \right),\\
s_3:(*) \rightarrow &\left(x,y,z,w-\frac{\alpha_3}{z},t;\alpha_0,\alpha_1,\alpha_2+\alpha_3,-\alpha_3 \right),\\
\pi:(*) \rightarrow &\left(1-x,-y,\sqrt{-1}z,\frac{w}{\sqrt{-1}},-t;\alpha_1,\alpha_0,\alpha_2,\alpha_3 \right).
\end{split}
\end{align}
\end{theorem}
These B{\"a}cklund transformations have Lie theoretic origin, similarity reduction of a Drinfeld-Sokolov hierarchy admits such a B{\"a}cklund symmetry.

\begin{proposition}\label{pro:3.4}
This system has the following invariant divisors\rm{:\rm}
\begin{center}
\begin{tabular}{|c|c|c|} \hline
parameter's relation & $f_i$ \\ \hline
$\alpha_0=0$ & $f_0:=x-1$  \\ \hline
$\alpha_1=0$ & $f_1:=x$  \\ \hline
$\alpha_2=0$ & $f_2:=yw^2-ty+1$  \\ \hline
$\alpha_3=0$ & $f_3:=z$  \\ \hline
\end{tabular}
\end{center}
\end{proposition}
We note that the system \eqref{31} admits a Riccati extension of the fifth Painlev\'e system as its particular solutions when $z=0$ with the parameter's relation $\alpha_3=0$.

\begin{proposition}
Let us define the following translation operators{\rm : \rm}
\begin{align}
\begin{split}
&T_1:=\pi s_0 s_2 s_3 s_2 s_0, \quad T_2:=\pi s_0 s_1 s_2 s_3 s_2, \quad T_3:=s_2 T_2 s_2.
\end{split}
\end{align}
These translation operators act on parameters $\alpha_i$ as follows$:$
\begin{align}
\begin{split}
T_1(\alpha_0,\alpha_1,\alpha_2,\alpha_3)=&(\alpha_0,\alpha_1,\alpha_2,\alpha_3)+(-1,1,0,0),\\
T_2(\alpha_0,\alpha_1,\alpha_2,\alpha_3)=&(\alpha_0,\alpha_1,\alpha_2,\alpha_3)+(1,1,-1,0),\\
T_3(\alpha_0,\alpha_1,\alpha_2,\alpha_3)=&(\alpha_0,\alpha_1,\alpha_2,\alpha_3)+(0,0,1,-2).
\end{split}
\end{align}
\end{proposition}

\begin{theorem}\label{th:3.2}
Let us consider a polynomial Hamiltonian system with Hamiltonian $K \in {\Bbb C}(t)[x,y,z,w]$. We assume that

$(B1)$ $deg(K)=6$ with respect to $x,y,z,w$.

$(B2)$ This system becomes again a polynomial Hamiltonian system in each coordinate system $r_i \ (i=0,1,3)${\rm : \rm}
\begin{align}
\begin{split}
r_0:&x_0=-((x-1)y-\alpha_0)y, \ y_0=\frac{1}{y}, \ z_0=z, \ w_0=w, \\
r_1:&x_1=-(xy-\alpha_1)y, \ y_1=\frac{1}{y}, \ z_1=z, \ w_1=w, \\
r_3:&x_3=x, \ y_3=y, \ z_3=-(zw-\alpha_3)w, \ w_3=\frac{1}{w}.
\end{split}
\end{align}
$(B3)$ In addition to the assumption $(B2)$, the Hamiltonian system in the coordinate $r_1$ becomes again a polynomial Hamiltonian system in the coordinate system $r_2${\rm : \rm}
\begin{equation*}
r_2:x_2=\frac{1}{x_1}, \ y_2=-\left((y_1+w_1^2-t)x_1+\alpha_2 \right)x_1, \ z_2=z_1-2x_1w_1, \ w_2=w_1.
\end{equation*}
Then such a system coincides with the system \eqref{31} with the polynomial Hamiltonian \eqref{32}.
\end{theorem}
By this theorem, we can also recover the parameter's relation \eqref{33}.

We note that the conditions $(B2)$ and $(B3)$ should be read that
\begin{align*}
&r_j(K) \quad (j=0,1,3), \quad r_2(r_1(K)+x_1)
\end{align*}
are polynomials with respect to $x,y,z,w$ or $x_1,y_1,z_1,w_1$.

\section{$D_4^{(3)}$ system}

In this section, we study a 2-parameter family of ordinary differential systems in dimension four with affine Weyl group symmetry of type $D_4^{(3)}$ given by
\begin{align}\label{41}
\begin{split}
\frac{dx}{dt}&=\frac{\partial H}{\partial y}, \quad \frac{dy}{dt}=-\frac{\partial H}{\partial x}, \quad \frac{dz}{dt}=\frac{\partial H}{\partial w}, \quad \frac{dw}{dt}=-\frac{\partial H}{\partial z}
\end{split}
\end{align}
with the polynomial Hamiltonian
\begin{align}\label{42}
\begin{split}
&(4t+3)(16t^2-12t+9)H=\\
&-12tx^2y^4+12(x+2\alpha_0 t)xy^3+3\{(\alpha_0+14\alpha_1+9\alpha_2)x-4\alpha_0^2 t\}y^2\\
&+\{-24tx^2-32t^2 x-3\alpha_0(5\alpha_0+14\alpha_1+9\alpha_2)\}y+24x\{x+(4\alpha_0+7\alpha_1+4\alpha_2)t\}\\
&-12tz^2w^4-12\{z-2(\alpha_0+\alpha_2)t\}zw^3+3\{8t^2 z^2-(\alpha_0+6\alpha_1-3\alpha_2)z-4(\alpha_0+\alpha_2)^2 t\}w^2\\
&+\{12tz^2-8(7\alpha_0+8\alpha_1+7\alpha_2)t^2 z+3(5-6\alpha_1-8\alpha_1^2-4\alpha_2)\}w\\
&-3\{(4t^3+1)z+(19\alpha_0+10\alpha_1+7\alpha_2)t\}z\\
&-12x\{-4txwy^3+(4xw+4txw^2-2tzw^2+zw+2t^2 z)y^2\\
&+(-2xw^2-7tz-4t^2 zw+zw^2+4tzw^3)y+(8tw+5)z\}+24\alpha_0^2 tyw\\
&+12\alpha_1 xw(2w-7y)+6\alpha_2 xyw\{4t(2w-y)-3\}\\
&+6\alpha_0 y\{4t(2x-z)w^2+(x-12t xy+2z)w+4t^2 z\}+24\alpha_0 \alpha_2 tyw.
\end{split}
\end{align}
Here $x,y,z$ and $w$ denote unknown complex variables, and $\alpha_0,\alpha_1,\alpha_2$ are complex parameters satisfying the relation:
\begin{equation}\label{43}
\alpha_0+2\alpha_1+\alpha_2=1.
\end{equation}
This is the first example which gave higher order Painlev\'e type systems of type $D_{4}^{(3)}$.

We remark that for this system we tried to seek its first integrals of polynomial type with respect to $x,y,z,w$. However, we can not find. Of course, the Hamiltonian $H$ is not the first integral.

This system can be obtained by connecting the invariant divisors $x-1$ and $z$ for the canonical variables $(x,y,z,w)$.

It is still an open question whether we make an explicit description of a confluence process from the system \eqref{31} to this system.

\begin{figure}[ht]
\unitlength 0.1in
\begin{picture}(34.95,8.70)(21.60,-16.60)
%
\special{pn 20}%
\special{ar 2585 1215 425 425  0.0000000 6.2831853}%
%
\special{pn 20}%
\special{ar 3960 1225 500 425  0.0000000 6.2831853}%
%
\special{pn 20}%
\special{ar 5230 1235 425 425  0.0000000 6.2831853}%
\put(24.1500,-13.0500){\makebox(0,0)[lb]{$x$}}%
\put(49.9000,-13.4000){\makebox(0,0)[lb]{$x+z$}}%
\put(35.1000,-13.3000){\makebox(0,0)[lb]{$yw^2-ty+1$}}%
%
\special{pn 20}%
\special{pa 3020 1240}%
\special{pa 3440 1240}%
\special{fp}%
%
\special{pn 20}%
\special{pa 4840 1040}%
\special{pa 4440 1040}%
\special{fp}%
\special{sh 1}%
\special{pa 4440 1040}%
\special{pa 4507 1060}%
\special{pa 4493 1040}%
\special{pa 4507 1020}%
\special{pa 4440 1040}%
\special{fp}%
%
\special{pn 20}%
\special{pa 4790 1260}%
\special{pa 4480 1260}%
\special{fp}%
\special{sh 1}%
\special{pa 4480 1260}%
\special{pa 4547 1280}%
\special{pa 4533 1260}%
\special{pa 4547 1240}%
\special{pa 4480 1260}%
\special{fp}%
%
\special{pn 20}%
\special{pa 4850 1460}%
\special{pa 4400 1460}%
\special{fp}%
\special{sh 1}%
\special{pa 4400 1460}%
\special{pa 4467 1480}%
\special{pa 4453 1460}%
\special{pa 4467 1440}%
\special{pa 4400 1460}%
\special{fp}%
\end{picture}
\label{fig:A63}
\caption{This figure denotes Dynkin diagram of type $D_4^{(3)}$. The symbol in each circle denotes the invariant divisors of the system \eqref{41} (see Theorem \ref{th:4.1}).}
\end{figure}

\begin{theorem}\label{th:4.1}
The system \eqref{41} admits the affine Weyl group symmetry of type $D_4^{(3)}$ as the group of its B{\"a}cklund transformations, whose generators $s_0,s_1,s_2$ defined as follows$:$ with {\it the notation} $(*):=(x,y,z,w,t;\alpha_0,\alpha_1,\alpha_2)$\rm{: \rm}
\begin{align}
\begin{split}
s_0:(*) \rightarrow &\left(x,y-\frac{\alpha_0}{x},z,w,t;-\alpha_0,\alpha_1+\alpha_0,\alpha_2 \right),\\
s_1:(*) \rightarrow &\left(x+\frac{\alpha_1(w^2-t)}{yw^2-ty+1},y,z+\frac{2\alpha_1 yw}{yw^2-ty+1},w,t;\alpha_0+\alpha_1,-\alpha_1,\alpha_2+3\alpha_1 \right),\\
s_2:(*) \rightarrow &\left(x,y-\frac{\alpha_2}{x+z},z,w-\frac{\alpha_2}{x+z},t;\alpha_0,\alpha_1+\alpha_2,-\alpha_2 \right).
\end{split}
\end{align}
\end{theorem}
These B{\"a}cklund transformations have Lie theoretic origin, similarity reduction of a Drinfeld-Sokolov hierarchy admits such a B{\"a}cklund symmetry.

\begin{proposition}
This system has the following invariant divisors\rm{:\rm}
\begin{center}
\begin{tabular}{|c|c|c|} \hline
parameter's relation & $f_i$ \\ \hline
$\alpha_0=0$ & $f_0:=x$  \\ \hline
$\alpha_1=0$ & $f_1:=yw^2-ty+1$  \\ \hline
$\alpha_2=0$ & $f_2:=x+z$  \\ \hline
\end{tabular}
\end{center}
\end{proposition}

\begin{theorem}\label{th:4.2}
Let us consider a polynomial Hamiltonian system with Hamiltonian $K \in {\Bbb C}(t)[x,y,z,w]$. We assume that

$(C1)$ $deg(K)=6$ with respect to $x,y,z,w$.

$(C2)$ This system becomes again a polynomial Hamiltonian system in each coordinate system $r_i \ (i=0,2)${\rm : \rm}
\begin{align}
\begin{split}
r_0:&x_0=-(xy-\alpha_0)y, \ y_0=\frac{1}{y}, \ z_0=z, \ w_0=w, \\
r_2:&x_2=x, \ y_2=y-w, \ z_2=-((z+x)w-\alpha_2)w, \ w_2=\frac{1}{w}.
\end{split}
\end{align}
$(C3)$ In addition to the assumption $(C2)$, the Hamiltonian system in the coordinate $r_1$ becomes again a polynomial Hamiltonian system in the coordinate system $r_0${\rm : \rm}
\begin{equation*}
r_1:x_1=\frac{1}{x_0}, \ y_1=-\left((y_0+w_0^2-t)x_0+\alpha_1 \right)x_0, \ z_1=z_0-2x_0w_0, \ w_1=w_0.
\end{equation*}
Then such a system coincides with the system \eqref{41} with the polynomial Hamiltonian \eqref{42}.
\end{theorem}
By this theorem, we can also recover the parameter's relation \eqref{43}.

We note that the conditions $(C2)$ and $(C3)$ should be read that
\begin{align*}
&r_j(K) \quad (j=0,2), \quad r_2(r_0(K)+x_0)
\end{align*}
are polynomials with respect to $x,y,z,w$ or $x_0,y_0,z_0,w_0$.

\section{Appendix}

In this appendix, we present a 2-parameter family of ordinary differential systems with affine Weyl group symmetry of type $G_2^{(1)}$ given by
\begin{align}\label{71}
\begin{split}
\frac{dx}{dt}&=\frac{\partial H}{\partial y}, \quad \frac{dy}{dt}=-\frac{\partial H}{\partial x}, \quad \frac{dz}{dt}=\frac{\partial H}{\partial w}, \quad \frac{dw}{dt}=-\frac{\partial H}{\partial z}
\end{split}
\end{align}
with the polynomial Hamiltonian
\begin{align}\label{72}
\begin{split}
H =&-3(5\alpha_0+18\alpha_1+27\alpha_2)txy^3-\frac{xy}{2t}+6\alpha_1(\alpha_0+6\alpha_1+9\alpha_2)ty^2\\
&-12t^2 z^2w^2-\left\{24(\alpha_1+\alpha_2)t^2-\frac{3}{2t} \right\}zw\\
&+12t^2 y^3zw^2-24txy^3zw+12x^2y^3z+24txz^2w-12x^2z^2-12\alpha_0 t^2y^3w\\
&+12(2\alpha_0+4\alpha_1+9\alpha_2)xy^2z+3(19\alpha_0+38\alpha_1+45\alpha_2)ty^2zw\\
&-6\alpha_1(\alpha_0+4\alpha_1+9\alpha_2)yz-4(\alpha_0-4\alpha_1-3\alpha_2)txz.
\end{split}
\end{align}
Here $x,y,z$ and $w$ denote unknown complex variables, and $\alpha_0,\alpha_1,\alpha_2$ are complex parameters satisfying the relation:
\begin{equation}\label{73}
\alpha_0+2\alpha_1+3\alpha_2=0.
\end{equation}
This is the first example which gave higher order Painlev\'e type systems of type $G_{2}^{(1)}$.

We remark that for this system we tried to seek its first integrals of polynomial type with respect to $x,y,z,w$. However, we can not find. Of course, the Hamiltonian $H$ and $z-ty$ are not its first integrals.

\begin{figure}[h]
\unitlength 0.1in
\begin{picture}(37.70,8.90)(7.40,-15.00)
%
\special{pn 20}%
\special{ar 1180 1050 440 440  0.0000000 6.2831853}%
%
\special{pn 20}%
\special{ar 2630 1060 440 440  0.0000000 6.2831853}%
%
\special{pn 20}%
\special{ar 4070 1060 440 440  0.0000000 6.2831853}%
%
\special{pn 20}%
\special{pa 1630 1050}%
\special{pa 2170 1050}%
\special{fp}%
%
\special{pn 20}%
\special{pa 3000 810}%
\special{pa 3660 810}%
\special{fp}%
\special{sh 1}%
\special{pa 3660 810}%
\special{pa 3593 790}%
\special{pa 3607 810}%
\special{pa 3593 830}%
\special{pa 3660 810}%
\special{fp}%
%
\special{pn 20}%
\special{pa 3080 1060}%
\special{pa 3590 1060}%
\special{fp}%
\special{sh 1}%
\special{pa 3590 1060}%
\special{pa 3523 1040}%
\special{pa 3537 1060}%
\special{pa 3523 1080}%
\special{pa 3590 1060}%
\special{fp}%
%
\special{pn 20}%
\special{pa 3000 1310}%
\special{pa 3670 1310}%
\special{fp}%
\special{sh 1}%
\special{pa 3670 1310}%
\special{pa 3603 1290}%
\special{pa 3617 1310}%
\special{pa 3603 1330}%
\special{pa 3670 1310}%
\special{fp}%
\put(9.8000,-11.7000){\makebox(0,0)[lb]{$z$}}%
\put(23.4000,-11.9000){\makebox(0,0)[lb]{$x-tw$}}%
\put(38.1000,-11.7000){\makebox(0,0)[lb]{$z-y^3$}}%
\end{picture}%
\label{fig:G21}
\caption{The symbol in each circle denotes the invariant divisors of this system (see Theorem \ref{th:7.1}).}
\end{figure}
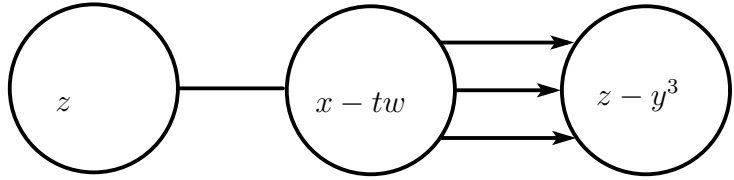

\begin{theorem}\label{th:7.1}
The system \eqref{71} admits the affine Weyl group symmetry of type $G_2^{(1)}$ as the group of its B{\"a}cklund transformations, whose generators $s_0,s_1,s_2$ defined as follows$:$ with {\it the notation} $(*):=(x,y,z,w,t;\alpha_0,\alpha_1,\alpha_2)$\rm{: \rm}
\begin{align}
\begin{split}
s_0:(*) \rightarrow &\left(x,y,z,w-\frac{\alpha_0}{z},t;-\alpha_0,\alpha_1+\alpha_0,\alpha_2 \right),\\
s_1:(*) \rightarrow &\left(x,y-\frac{\alpha_1}{x-tw},z-\frac{\alpha_1 t}{x-tw},w,t;\alpha_0+\alpha_1,-\alpha_1,\alpha_2+\alpha_1 \right),\\
s_2:(*) \rightarrow &\left(x-\frac{3 \alpha_2 y^2}{z-y^3},y,z,w-\frac{\alpha_2}{z-y^3},t;\alpha_0,\alpha_1+3\alpha_2,-\alpha_2  \right).
\end{split}
\end{align}
\end{theorem}
Since these B{\"a}cklund transformations have Lie theoretic origin, similarity reduction of a Drinfeld-Sokolov hierarchy admits such a B{\"a}cklund symmetry.

\begin{proposition}
This system has the following invariant divisors\rm{:\rm}
\begin{center}
\begin{tabular}{|c|c|c|} \hline
parameter's relation & $f_i$ \\ \hline
$\alpha_0=0$ & $f_0:=z$  \\ \hline
$\alpha_1=0$ & $f_1:=x-tw$  \\ \hline
$\alpha_2=0$ & $f_2:=z-y^3$  \\ \hline
\end{tabular}
\end{center}
\end{proposition}

\begin{theorem}\label{pro:7.2}
Let us consider a polynomial Hamiltonian system with Hamiltonian $K \in {\Bbb C}(t)[x,y,z,w]$. We assume that

$(F1)$ $deg(K)=6$ with respect to $x,y,z,w$.

$(F2)$ This system becomes again a polynomial Hamiltonian system in each coordinate system $r_i \ (i=0,1,2)${\rm : \rm}
\begin{align}
\begin{split}
r_0:&x_0=x, \ y_0=y, \ z_0=-(zw-\alpha_0)w, \quad w_0=\frac{1}{w},\\
r_1:&x_1=-((x-tw)y-\alpha_1)y, \ y_1=\frac{1}{y}, \ z_1=z-ty, \ w_1=w,\\
r_2:&x_2=x-3y^2 w, \ y_2=y, \ z_2=-((z-y^3)w-\alpha_2)w, \ w_2=\frac{1}{w}.
\end{split}
\end{align}
Then such a system coincides with the system
\begin{align}
\begin{split}
\frac{dx}{dt}&=\frac{\partial K}{\partial y}, \quad \frac{dy}{dt}=-\frac{\partial K}{\partial x}, \quad \frac{dz}{dt}=\frac{\partial K}{\partial w}, \quad \frac{dw}{dt}=-\frac{\partial K}{\partial z}
\end{split}
\end{align}
with the polynomial Hamiltonian
\begin{align}\label{75}
\begin{split}
K =&H+\sum_{i=1}^6 a_i(z-ty)^i \quad (a_i \in {\Bbb C}(t)).
\end{split}
\end{align}
\end{theorem}
We note that the condition $(F2)$ should be read that
\begin{align*}
&r_j(K) \quad (j=0,2), \quad r_1(K-yw)
\end{align*}
are polynomials with respect to $x,y,z,w$.

\end{document}